\documentclass[12pt, reqno]{amsart}
\usepackage[margin=1.3in]{geometry}

\usepackage{bbm,amsmath,amssymb,amsfonts}
\usepackage[comma, sort&compress]{natbib}
\usepackage[pdfpagemode=UseNone,bookmarksopen=false,colorlinks=true]{hyperref}
\usepackage{pstricks}
\usepackage{enumerate}
\usepackage{tikz}
\usetikzlibrary{matrix}
\usepackage{bbm}

\hypersetup{ colorlinks=true, linkcolor=blue, citecolor=blue,
	filecolor=blue, urlcolor=blue}
\usepackage[utf8]{inputenc}
\usepackage[english]{babel}
\newcommand{\zd}{\mathbb{Z}^d}
\newcommand{\td}{\mathbb{T}_n^d}
\newcommand{\ds}{\displaystyle}
\DeclareMathOperator{\prob}{\mathbf{P}}

\newtheorem{theorem}{Theorem}
\numberwithin{theorem}{section} 
\newtheorem{lemma}[theorem]{Lemma}

\newtheorem{corollary}[theorem]{Corollary}

\newtheorem{definition}[theorem]{Definition}

\usepackage{bbm}

\newcommand{\Z}{\mathbb Z}

\begin{document}
	
	\title[Extremes of Zero-Average DGFF]{Extremal process of the zero-average Gaussian Free Field for $d\ge 3$}
	\author[S.~Das]{Sayan Das}
	\address{Columbia University, 2990 Broadway, New York, NY 10027}
	\email{sayan.das@columbia.edu}
	\author[R.~S.~Hazra]{Rajat Subhra Hazra}
	\address{Stat-Math Unit, Indian Statistical Institute, Kolkata}
	\email{rajatmaths@gmail.com}
	\date{}

	\begin{abstract}
		We consider the Gaussian free field on the torus whose covariance kernel is given by the zero-average Green's function. We show that for dimension $d\ge 3$, the extremal point process associated with this field converges weakly to a Poisson random measure. As an immediate corollary the maxima of the field converges after appropriate centering and scaling to the Gumbel distribution.
	\end{abstract}
	
	
	\maketitle
	\section{\textbf{Introduction}}
	
	The Gaussian free field is an important example of a random interface.  In $\Z^d$ the Gaussian free field on a box of side length $n$ is defined as a centered Gaussian process whose covariance is given by the Green's function of the simple random walk on $\Z^d$ (conditioned to be killed upon exiting the box). In $d=2$, it falls under the rich class of log-correlated models and due to its close connections with the branching random walk, it has been an important object of study (\cite{ding:roy:ofer}). In $d\ge 3$, a deep structure appears when one considers the level set percolation (\cite{PFASS}). The continuum Gaussian free field, especially in $d=2$, also plays a crucial role in SLE theory, due to the natural conformal invariance property (\cite{schramm2009contour}). In this short note, we are interested in the study of extremes of a Gaussian free field on a torus. In $\Z^d$ the picture of the extremal process has become clearer due to important contributions in the works of \cite{biskup2016full, BZ2012, bramson-zeitouni, ding:roy:ofer, chiarini2015note, chiarini2016extremes}. For a comprehensive review we refer to \cite{biskup2017extrema}.
	
	The zero-average Gaussian free field is a centered Gaussian field, indexed by a discrete torus whose covariance is given by the zero average Green's function (see \eqref{eq:zeroav:green} for the precise definition). The zero-average Gaussian fields are known to be related to occupation measures of charged particles (see \cite[Chapter 14.6.2]{aldous2000reversible} for more details). The nomenclature of ``zero average" was introduced recently by \cite{angelo}, in the context of  the level set percolation. Due to the lack of boundary on the torus, there are other approaches to define the Green's function on the torus which leads to  different kinds of Gaussian free field. In $d=2$ they were studied in ~\cite{chatterjee2016superconcentration}, and~\cite{BZ2012}. 
	
	The zero average Green's function is closely related to the Green's function on $\Z^d$. In \cite{angelo} a nice coupling between the fields on the discrete torus and $\Z^d$ was given. We show that this coupling gives us an easy and natural way to derive the scaling limit of the extremal process. We show the point process of the Gaussian free field converges to a Poisson random measure. The behaviour of the point process turns out to be similar to the case of $\Z^d$ in $d\ge 3$ and hence like a collection of  i.i.d.~Gaussians.  As a corollary, the weak convergence of the maxima (after centering and scaling) to the standard Gumbel distribution follows. The long-range correlation does not affect the extremal process. At this moment, we cannot confirm whether this can be extended to the more interesting $d=2$ case. It would be interesting to see if a Cox-cluster process (similar to \cite{biskup2016full}) appears there also.
	
	The outline of the article is as follows. Since the main ingredient of the proof is a coupling result by \cite{angelo}, we recall the notations and some of the important results from there in the next section. In the section~\ref{4sec2}, we state the main idea of the proof and then finally in section~\ref{4sec3} we give the complete arguments for the intermediate steps.

	\section{Notations and main result}
	We borrow the notations and setup from \cite{angelo} with a slight modification. We consider the discrete tori $\td := (\mathbb{Z}/n\mathbb{Z})^d$, $n \ge 1$, and $\mathbb{Z}^d$ for $d \ge 3$ endowed with the usual graph structure. 
	We let $\Pi_n : \zd \to \td$ denote the canonical projection. If $x \in \td$, we write $\widehat{x} \in \zd$ for the unique element
	of $\Pi_n^{-1} (\{x\})\cap \left[0,n-1\right]^d$. If $U \subset \td$, we similarly write $\widehat{U} := \{\widehat{x} \in \zd \mid x \in U\} \subset \zd$. Note
	that the map $x \to \widehat{x}$ from $\td$ to $\left[0,n-1\right]^d \cap \zd$ is a bijection with inverse $\Pi_n\mid_{\left[0,n-1\right]^d \cap \zd}$.
	
	We consider the simple random walks on $\td$ and $\zd$. We write $P_x^{\td}$ and $P_x^{\zd}$ to denote the law of the simple random walks on $\td$ and $\zd$ starting at $x$ respectively. Moreover, we shall write $E_x^{\td}$ and $E_x^{\zd}$ for the corresponding expectations respectively. The canonical process for both discrete-time walks is denoted by $(X_k)_{k\ge0}$. For the continuous-time walks with i.i.d.~$\operatorname{Exp}(1)$ holding times, we write $(\bar{X}_t)_{t\ge0}$.

	The Green's function $g_{\zd}(\cdot,\cdot)$ for the simple random walk on $\zd$ is
	$$g_{\zd}(x,y):=E_x^{\zd}\left[\sum_{k=0}^{\infty} \mathbb{I}_{\{X_k=y\}}\right] = \sum_{k=0}^{\infty} P_x^{\zd}[X_k=y] \ \ \ \mbox{ for }x,y\in \zd.$$
	
	Note that the above is finite as we deal with the case $d\ge 3$. Additionally, it is symmetric, positive,  and satisfies $g_{\zd}(x,y)=g_{\zd}(x-y,0)$. 
	
	The zero-average Green's function $G_{\td}(\cdot,\cdot)$ associated with the simple random walk on $\td$ is given by
	\begin{equation}\label{eq:zeroav:green}
	G_{\td}(x,y):=\int_0^{\infty} \left(P_x^{\td}[\bar{X}_t=y]-\frac1{n^d} \right)\,dt \ \ \ \mbox{ for }x,y\in \td.
	\end{equation}
	It turns out that $G_{\td}(\cdot,\cdot)$ is also symmetric, finite and positive semi-definite and satisfies $G_{\td}(x,y)=G_{\td}(x-y,0)$. We now recall the definitions of Gaussian free field on $\td$ and $\Z^d$.
	
	\begin{definition}
		
		
		The zero-average Gaussian free field $(\Psi_{\td}(x))_{x\in\td}$ is a centered Gaussian field on $\td$ with covariance $$\mathbb{E}^{\td}[\Psi_{\td}(x)\Psi_{\td}(y)]= G_{\td}(x,y) \text{ for all $x, y \in \td$}.$$ 
		The law of $(\Psi_{\td}(x))_{x\in\td}$ on $\mathbb{R}^{\td}$ is denoted by $\mathbb{P}^{\td}$.
		
		On the other hand, on $\zd$ we have the infinite-volume Gaussian free field $(\varphi_{\zd}(x))_{x\in\zd}$ with law $\mathbb{P}^{\zd}$ on $\mathbb{R}^{\zd}$. It is a centered Gaussian field on $\zd$ with covariance structure given by $$\mathbb{E}^{\zd}[\varphi_{\zd}(x)\varphi_{\zd}(y)]= g_{\zd}(x,y) \text{ for all $x, y \in \zd$}.$$
	\end{definition}
	
	We remind the reader that the nomenclature of zero-average was introduced in \cite{angelo} due to the following property
	$$\textrm{Var}_{\mathbb{P}^{\td}}\left(\sum_{x\in \td} \Psi_{\td}(x)\right)=0.
	$$
	We now state the main results of this note on the extremal process and maxima of the zero-average Gaussian free field. 
	\subsection{Main results}
	We let $E=(\mathbb{R}/\mathbb{Z})^d\times (-\infty,\infty]$ and denote $\mathcal{M}_p(E)$ to be the space of all Radon point measures on $E$ endowed with the topology of vague convergence. We define the following sequence of extremal point processes on $E$ associated with zero-average Gaussian free field.
	\begin{equation} \label{pp1}
	\eta_n:=\sum_{\alpha\in \td} \delta_{\left(\frac{\alpha}{n},\frac{\Psi_{\mathbb{T}_n^d}(\alpha)-b_N}{a_N}\right)}(\cdot)
	\end{equation}
	where $\delta_x(\cdot), x \in E$, is the point measure that gives mass one to a set containing $x$ and zero otherwise, and
	\begin{equation}\label{con}
	b_N=\sqrt{g_{\zd}(0,0)}\left[\sqrt{2\log N}-\frac{\log \log N+\log 4\pi}{2\sqrt{2\log N}}\right], \ \ \ \ a_N = g_{\zd}(0,0)b_N^{-1}, \ \ \ N=n^d.
	\end{equation}
	
	Our main result is 
	
	\begin{theorem} \label{zagff}
		For the sequence of point processes $\eta_n$ defined in \eqref{pp1} we have that
		$$\eta_n\stackrel{d}{\to} \eta$$
		as $n \to \infty$, where $\eta$ is a Poisson random measure on $E$ with intensity measure given by $dt\otimes e^{-z} dz$	where $dt\otimes dz$ is the Lebesgue measure on $E$, and $\stackrel{d}{\to}$ is the convergence in distribution on the space $\mathcal{M}_p(E)$.
	\end{theorem}
	
	The proof is based  on the point process convergence result on $\zd$ of \cite{chiarini2015note} and a coupling lemma of \cite{angelo} that allow us to compare the Gaussian free field on $\zd$ and the zero-average Gaussian free field on $\td$. As a corollary, we obtain the limiting distribution of the maximum of the Gaussian free field on the torus.
	
	\begin{corollary}
		The maximum of the zero-average Gaussian free field on the torus belongs to the domain of attraction of Gumbel distribution. In particular, we have for all $z\in \mathbb{R}$
		\begin{equation}\lim_{n \to \infty} P\left(\max_{\alpha\in \td} \Psi_{\td}(\alpha) \le a_Nz+b_N \right)=\exp(-e^{-z}).\end{equation}
	\end{corollary}

	\subsection{Some known facts}	
	To keep the article self-contained, in this subsection we write down some known estimates of the Green's function on $\Z^d$ and $\td$ and also recall the main ingredient of the proof, namely the coupling result from \cite{angelo}. 	
	
	Let us write $|x|$ to denote the Euclidean norm of $x\in \Z^d$. In $d\ge 3$, $g_{\zd}(x,y)$ has a polynomial decay of order $|x-y|^{2-d}$ as $|x-y|\to\infty$ and the following lemma asserts this. The estimate will be crucial in our proofs. 
	
	\begin{lemma}[Theorem 1.5.4 of \cite{lawler2013intersections}] \label{two} For any $x,y \in \zd$, it holds that	$$c_d|x-y|^{2-d} \le g_{\zd}(x,y)\le C_d|x-y|^{2-d}$$
		where $0<c_d\le C_d <\infty$.	
	\end{lemma} 
	For $V \subset \zd$, the Green's function $g_{\zd}^V(\cdot,\cdot)$ of the simple random walk on $\zd$ killed when exiting $V$ is
	$$g_{\zd}^V(x,y):=E_x^{\zd}\left[\sum_{0\le k<T_V} \mathbb{I}_{\{X_k=y\}}\right]= \sum_{k=0}^{\infty} P_x^{\zd}[X_k=y,k<T_V] \ \ \ \mbox{ for }x,y\in \zd,$$
	where $T_V:=\inf\{k\ge 0 \mid X_k \notin V\}$ is the exit time from $V$. A key fact for the Green's function is the spatial Markov property which is stated in the next lemma. 
	\begin{lemma}[Proposition 4.6.2(a) of \cite{lawler2010random}] \label{one} For $V\subsetneq  \zd$ and $x,y\in\zd$, we have
		$$g_{\zd}(x,y)=g_{\zd}^{V}(x,y)+E_x^{\zd}[g_{\zd}(X_{T_V},y)\mathbb{I}_{\{T_V<\infty\}}]. $$
	\end{lemma}

	We define for $U \subsetneq \td$ the Green's function of the simple random walk on $\td$ killed when exiting $U$, which is
	$$g_{\td}^U(x,y):=E_x^{\td}\left[\sum_{0\le k<T_U} \mathbb{I}_{\{X_k=y\}}\right]= \sum_{k=0}^{\infty} P_x^{\td}[X_k=y,k<T_U] \ \ \ \mbox{ for }x,y\in \td.$$
	We now state some properties of the Green's function $G_{\td}(\cdot,\cdot)$. 
	\begin{lemma}[Lemma 1.3 of \cite{angelo}] \label{thr} Assume $U\subsetneq \td$. Then it holds that
		$$G_{\td}(x,y)=g_{\td}^U(x,y)+E_x^{\td}[G_{\td}(X_{T_U},y)]-\frac1{n^d}E_x^{\td}[T_U] \ \ \mbox{for all }x,y\in\td.$$
	\end{lemma}
	
	Let us denote the usual graph distance on $\td$ by $d_{\td}(\cdot,\cdot)$. 
	\begin{lemma}[Proposition 1.4 of \cite{angelo}] \label{four} For all $n\ge 1$ and $x,y\in \td$ it holds that
		$$|G_{\td}(x,y)|\le c(\log (n))^{\frac{3d}2}d_{\td}(x,y)^{2-d}+c^\prime n(n\log(n))^{d+1}e^{-c^{\prime\prime}(\log(n))^2}$$
		where $c, c^\prime$ and $c^{\prime\prime}$ are constants depending on $d$ only.
	\end{lemma}
	Lemma \ref{thr} is the zero-average Green's function analog of Lemma \ref{one}. Lemma \ref{four} shows that the Green's function $G_{\td}(x,y)$ goes to zero when the points $x,y$ are far apart in $d_{\td}$ as $n\to\infty$.

	We close this section with a powerful coupling result.
	
	\begin{lemma}[Theorem 2.3 of \cite{angelo}] \label{angl}
		Let $R_n=(n^{3/4},n-n^{3/4}]^d \cap \zd$. For any $n\ge 1$ there exists a coupling $\mathbb{Q}_n$ of $\Psi_{\td}$ and $\varphi_{\zd}$ such that for all $\epsilon>0$
		$$\mathbb{Q}_n\left[\sup_{x\in R_n}\left|\Psi_{\td}(\Pi_n(x))-\varphi_{\zd}(x)\right|\ge \epsilon\right]\le 4n^d\exp(-c_1\epsilon^2n^{c_2})$$
		where $c_1,c_2$ are positive constants.
	\end{lemma} 
	The above result is true for sets of the form $R_n=(n^{\beta},n-n^{\beta}]^d \cap \zd$ where $\beta \in (\frac12,1)$. However, for our proof it suffices to consider $\beta=3/4$.

	
	\section{\bf Outline of Proof of Theorem \ref{zagff}}\label{4sec2}
	We consider the space $\widetilde{E}:=[0,1]^d\times (-\infty,\infty]$. Let $V_n=[0,n-1]^d\cap \zd$ be the $d$-dimensional box of side length $n$ on $\zd$. Let $C_c^{+}(\widetilde{E})$ be the collection of all non-negative continuous real-valued functions on $\widetilde{E}$ with compact support. We denote the space of all Radon point measures on $\widetilde{E}$ as $\mathcal{M}_p(\widetilde{E})$ endowed with the topology of vague convergence. This topology is known to be metrizable by the metric $$\widetilde\rho(\mu,\nu):= \sum_{i=1}^{\infty} 2^{-i}\min (|\mu(h_i)-\nu(h_i)|, 1),$$ where
	$\{h_i\}_{i\ge1}$ is a suitably chosen subset of $C_c^{+}(E)$ consisting only of Lipschitz functions (see Proposition 3.17 and Lemma 3.11 of \cite{resnick:1987}). We consider the following two point processes as random elements in $\mathcal{M}_p(\widetilde{E})$ :
	$$\widetilde\eta_n:=\sum_{\alpha\in V_n} \delta_{\left(\frac{\alpha}{n},\frac{\Psi_{\mathbb{T}_n^d}(\Pi_n(\alpha))-b_N}{a_N}\right)}(\cdot), \ \ \ \ \ \zeta_n:=\sum_{\alpha\in V_n} \delta_{\left(\frac{\alpha}{n},\frac{\varphi_{\zd}(\alpha)-b_N}{a_N}\right)}(\cdot).$$
	
	It is known from \cite{chiarini2015note} that the point process $\zeta_n$ converges weakly in $\mathcal{M}_p(\widetilde{E})$. In fact,
	\begin{equation}\label{rsh} \zeta_n\stackrel{d}{\to} \widetilde\eta\end{equation}
	as $n \to \infty$, where $\widetilde\eta$ is a Poisson random measure on $\widetilde{E}$ with intensity measure given by $dt\otimes e^{-z} dz$	where $dt\otimes dz$ is the Lebesgue measure on $\widetilde{E}$.
	
	Our proof relies on the following lemma which essentially establish that the asymptotic behaviour of $\widetilde{\eta}_n$ is same as that of $\zeta_n$. 
	
	\begin{lemma}\label{slut} There exist a coupling $\prob_n$ of $\Psi_{\td}$ and $\varphi_{\zd}$ such that for all $\epsilon>0$ we have
		$$\lim_{n\to \infty} \prob_n \left(\widetilde\rho\left(\widetilde\eta_n,\zeta_n\right)\ge \epsilon\right)=0.$$
	\end{lemma}
	
	We postpone the proof of the lemma to next section. We now complete the proof of Theorem \ref{zagff} using Lemma \ref{slut}. Towards this end, we define $\Pi: [0,1]^d \to (\mathbb{R}/\mathbb{Z})^d$ to be the natural projection map. We define a function $h: \mathcal{M}_p(\widetilde{E})\to \mathcal{M}_p({E})$ as follows. For any $\widetilde\mu=\sum_i \delta_{(\alpha_i,x_i)}\in \mathcal{M}_p(\widetilde{E})$, define $h(\widetilde\mu):= \mu$ where $\mu=\sum_i \delta_{(\Pi(\alpha_i),x_i)}\in \mathcal{M}_p(E)$. Note that $h$ is continuous. To see this, let $\widetilde\mu=\sum_i \delta_{(\alpha_i,x_i)}\in \mathcal{M}_p(\widetilde{E})$. Take a sequence $\widetilde\mu_n=\sum_i \delta_{(\alpha_{i,n},x_{i,n})}\in \mathcal{M}_p(\widetilde{E})$ converging to $\widetilde\mu$ with respect to the metric $\widetilde\rho$. We will show that $h(\widetilde\mu_n)$ converges to $h(\widetilde\mu)$ with respect to the vague topology on $\mathcal{M}_p({E})$. Take $f:E\to \mathbb{R}$ to be a compactly supported continuous function. Define $g : \widetilde{E} \to \mathbb{R}$ as $g(\alpha,x)= f(\Pi(\alpha),x)$. 
	Clearly, $g$ is also a compactly supported continuous function. Hence
	$$h(\widetilde\mu_n)(f)=\sum_i f(\Pi(\alpha_{i,n}),x_{i,n}) = \sum_i g(\alpha_{i,n},x_{i,n})=\widetilde\mu_n(g) \to \widetilde\mu(g)=h(\widetilde\mu)(f).$$
	This establishes the continuity of $h$. Observe that by Slutsky's Theorem, one can combine Lemma \ref{slut} and \eqref{rsh}, to get $\widetilde\eta_n \stackrel{d}{\to} \widetilde\eta$. Note that $h(\widetilde{\eta}_n)=\eta_n$ and $h(\widetilde{\eta})$ has the same distribution as $\eta$. Hence by Continuous Mapping Theorem we have
	$$\eta_n \stackrel{d}{\to} \eta.$$
	This completes the proof of Theorem \ref{zagff} modulo Lemma \ref{slut}.
	
	\section{\bf Proof of Lemma \ref{slut}}\label{4sec3}
	We consider the coupling $\mathbb{Q}_n$ of 	$\Psi_{\td}$ and $\varphi_{\zd}$ as described in Lemma \ref{angl}. In fact, we will show that this coupling $\mathbb{Q}_n$ is indeed the coupling $\mathbf{P}_n$ of Lemma \ref{slut}. For simplicity, let us write $$\psi(\alpha):=\dfrac{\Psi_{\mathbb{T}_n^d}(\Pi_n(\alpha))-b_N}{a_N}, \mbox{ and } \phi(\alpha):=\dfrac{\varphi_{\zd}(\alpha)-b_N}{a_N}$$
	where $N=n^d$ is the size of the torus $\td$. With these notations in hand, we see that $$\ds \widetilde\eta_n=\sum_{\alpha\in V_n} \delta_{\left(\frac{\alpha}{n},\psi(\alpha)\right)}(\cdot), \mbox{ and }\ds \zeta_n=\sum_{\alpha\in V_n} \delta_{\left(\frac{\alpha}{n},\phi(\alpha)\right)}(\cdot).$$ 
	
	Take a Lipschitz function $g\in C_c^{+}(\widetilde{E})$. Assume that the support of $g$ is contained in $[0,1]^d \times  (\delta,\infty]$ for some $\delta \in \mathbb{R}$. By the definition of the vague metric $\widetilde{\rho}$, it suffices to prove that, for every $\epsilon>0$,
	\begin{equation} \label{geqn}
	\limsup_{ n \to \infty} \mathbb{Q}_n(|\widetilde{\eta}_n(g)-\zeta_n(g)|\ge \epsilon) = 0.
	\end{equation}
	Let $R_n=(n^{3/4},n-n^{3/4}]^d \cap \zd$.  The main idea is to show that the contributions essentially come when the field is restricted to $R_n$, where we can also apply Lemma \ref{angl}. Towards this end, we define the following events
	$$\mathbf{A}_n= \{\psi(\alpha) \le \delta \mid \ \forall \ \alpha \in V_n\setminus R_n \}, \mbox{ and } \mathbf{B}_n= \{\phi(\alpha) \le \delta \mid \ \forall \ \alpha \in V_n\setminus R_n \}. $$
	Moreover, for each $0<\gamma<1$, we define
	$$\mathbf{C}_n(\gamma)= \{|\psi(\alpha)-\phi(\alpha)| \le \gamma \mid \ \forall \ \alpha \in R_n \}.$$ 
	The following lemma establishes that with high probability $\mathbf{A}_n, \mathbf{B}_n$, and $\mathbf{C}_n(\gamma)$ occur.
	\begin{lemma} \label{tech} For each $\gamma \in (0,1)$, we have
		$$\limsup_{ n \to \infty} \mathbb{Q}_n(\mathbf{A}_n^{\mathsf{c}})=\limsup_{ n \to \infty} \mathbb{Q}_n(\mathbf{B}_n^{\mathsf{c}})=\limsup_{ n \to \infty} \mathbb{Q}_n\left[(\mathbf{C}_n(\gamma))^{\mathsf{c}}\right]=0.$$
	\end{lemma}
	The proof of this technical result is postponed to the next section. Assuming Lemma \ref{tech}, the proof of Lemma \ref{slut} can be completed as follows. Note that
	\begin{equation} \label{ineq} \begin{aligned}\mathbb{Q}_n(|\widetilde{\eta}_n(g)-\zeta_n(g)|\ge \epsilon) & \le \mathbb{Q}_n(|\widetilde{\eta}_n(g)-\zeta_n(g)|\ge \epsilon,\mathbf{A}_n,\mathbf{B}_n,\mathbf{C}_n(\gamma))+\\ & \hspace{1.3cm} +\mathbb{Q}_n(\mathbf{A}_n^{\mathsf{c}})+\mathbb{Q}_n(\mathbf{B}_n^{\mathsf{c}})+\mathbb{Q}_n\left[(\mathbf{C}_n(\gamma))^{\mathsf{c}}\right].\end{aligned}\end{equation}
	Lemma \ref{tech} will imply that the last three terms of the right side of \eqref{ineq} are asymptotically zero. Observe that $g\left(\frac{\alpha}{n},\psi(\alpha)\right)$ is zero whenever $\psi(\alpha)\le \delta$. Hence we can write
	
	\begin{equation} \label{exp1}
	\widetilde{\eta}_n(g)=\sum_{\alpha\in  R_n} g\left(\frac{\alpha}{n},\psi(\alpha)\right)+\sum_{\alpha\in V_n\setminus R_n} g\left(\frac{\alpha}{n},\psi(\alpha)\right)\mathbbm{1}_{\mathbf{A}_n^{\mathsf{c}}}.
	\end{equation}
	Similarly one can use $\mathbf{B}_n$ to decompose $\zeta_n(g)$. We have
	
	\begin{equation} \label{exp2}
	\zeta_n(g)=\sum_{\alpha\in R_n} g\left(\frac{\alpha}{n},\phi(\alpha)\right)+\sum_{\alpha\in V_n\setminus R_n} g\left(\frac{\alpha}{n},\phi(\alpha)\right)\mathbbm{1}_{\mathbf{B}_n^{\mathsf{c}}}.
	\end{equation}
	The above decomposition is crucial in bounding the first term on the right hand side of \eqref{ineq}. Conditioning on the events $\mathbf{A}_n$ and $\mathbf{B}_n$, one can get rid of the second summand appearing in both \eqref{exp1} and \eqref{exp2}. Hence we have
	\begin{equation} \label{pr1}
	\begin{aligned}
	& \mathbb{Q}_n\left(\left|\widetilde{\eta}_n(g)-\zeta_n(g)\right|\ge \epsilon,\mathbf{A}_n,\mathbf{B}_n,\mathbf{C}_n(\gamma)\right)  \\ & \le \mathbb{Q}_n\left(\left|\sum_{\alpha\in R_n} g\left(\frac{\alpha}{n},\psi(\alpha)\right)-\sum_{\alpha\in R_n} g\left(\frac{\alpha}{n},\phi(\alpha)\right)\right|\ge \epsilon,\mathbf{C}_n(\gamma)\right) \\ &  \le \mathbb{Q}_n\left(\sum_{\alpha\in R_n}\left| g\left(\frac{\alpha}{n},\psi(\alpha)\right)- g\left(\frac{\alpha}{n},\phi(\alpha)\right)\right|\ge \epsilon,\mathbf{C}_n(\gamma)\right) \\ & \le \mathbb{Q}_n\left(\sum_{\alpha\in R_n}\left| g\left(\frac{\alpha}{n},\psi(\alpha)\right)- g\left(\frac{\alpha}{n},\phi(\alpha)\right)\right|\mathbbm{1}_{\mathbf{C}_n(\gamma)}\ge \epsilon\right). 
	\end{aligned}
	\end{equation}

	We next show that the random variable 
	\begin{equation}
	\label{rv1} \ds \sum_{\alpha\in R_n}\left| g\left(\frac{\alpha}{n},\psi(\alpha)\right)- g\left(\frac{\alpha}{n},\phi(\alpha)\right)\right|\mathbbm{1}_{\mathbf{C}_n(\gamma)}
	\end{equation}
	appearing in the last line of \eqref{pr1} can be bounded by a suitable random variable. Towards this end, we suppose $\phi(\alpha) \le \delta-\gamma$ and $|\psi(\alpha)-\phi(\alpha)| \le \gamma$. Then $\psi(\alpha) \le \delta$ and as a consequence, both $g\left(\frac{\alpha}{n},\psi(\alpha)\right)$ and $g\left(\frac{\alpha}{n},\phi(\alpha)\right)$ are zero. Furthermore, under our assumption, we have 
	$$\left|g\left(\frac{\alpha}{n},\psi(\alpha)\right)-g\left(\frac{\alpha}{n},\phi(\alpha)\right)\right| \le \|g\| |\psi(\alpha)-\phi(\alpha)| \le \|g\|\gamma$$
	where $\|g\|$ denotes the Lipschitz constant. These estimates allow us to give an upper bound on \eqref{rv1}. We have
	
	\begin{equation} \label{pr2}
	\begin{aligned}
	& \sum_{\alpha\in R_n}\left| g\left(\frac{\alpha}{n},\psi(\alpha)\right)- g\left(\frac{\alpha}{n},\phi(\alpha)\right)\right|\mathbbm{1}_{\mathbf{C}_n(\gamma)} \\ & = \sum_{\alpha\in R_n}\left| g\left(\frac{\alpha}{n},\psi(\alpha)\right)- g\left(\frac{\alpha}{n},\phi(\alpha)\right)\right|\mathbb{I}_{\phi(\alpha)>\delta-\gamma}\mathbbm{1}_{\mathbf{C}_n(\gamma)} \\ & \le \gamma\|g\|\sum_{\alpha\in R_n} \mathbb{I}_{\phi(\alpha)>\delta-\gamma}\mathbbm{1}_{\mathbf{C}_n(\gamma)}  \\ & \le \gamma\|g\|\zeta_n([0,1]^d\times (\delta-\gamma,\infty]) \\ & \le \gamma\|g\|\zeta_n([0,1]^d\times (\delta-1,\infty]).
	\end{aligned}
	\end{equation}
	
	Now using the point process result in $\eqref{rsh}$, we see that $\zeta_n([0,1]^d\times (\delta-1,\infty])$ converges weakly to  $\widetilde\eta([0,1]^d\times (\delta-1,\infty])$, which is  finite almost surely. Thus combining \eqref{pr1} and \eqref{pr2}, we get that
	\begin{equation}\label{pr3}
	\begin{aligned}
	& \limsup_{n \to \infty} \mathbb{Q}_n(|\widetilde{\eta}_n(g)-\zeta_n(g)|\ge \epsilon,\mathbf{A}_n,\mathbf{B}_n,\mathbf{C}_n(\gamma)) \\ & \le  \limsup_{ n \to \infty} \mathbb{Q}_n\left(\sum_{\alpha\in R_n}\left| g\left(\frac{\alpha}{n},\psi(\alpha)\right)- g\left(\frac{\alpha}{n},\phi(\alpha)\right)\right|\mathbbm{1}_{\mathbf{C}_n(\gamma)}\ge \epsilon\right)  \\ & \le \limsup_{ n \to \infty} \mathbb{Q}_n (\gamma\|g\|\zeta_n([0,1]^d\times (\delta-1,\infty])\ge\epsilon) \\ & = \limsup_{ n \to \infty} \mathbb{P}^{\zd} (\gamma\|g\|\zeta_n([0,1]^d\times (\delta-1,\infty])\ge\epsilon) \\ &\le \prob (\gamma\|g\| \widetilde{\eta}([0,1]^d\times (\delta-1,\infty])\ge\epsilon).
	\end{aligned}
	\end{equation}

	The last inequality follows from Portmanteau theorem. Here we assume $(\Omega,\mathcal{F},\prob)$ is the probability space where $\widetilde{\eta}$ is defined. Thus by taking limsup on both sides of \eqref{ineq} and applying Lemma \ref{tech}, and the bound in \eqref{pr3}, we get
	\begin{equation}\label{pr4}
	\limsup_{ n \to \infty}  \mathbb{Q}_n(|\widetilde{\eta}_n(g)-\zeta_n(g)|\ge \epsilon) \le \prob (\gamma\|g\| \widetilde{\eta}([0,1]^d\times (\delta-1,\infty])\ge\epsilon).
	\end{equation}
	Since $\gamma$ is arbitrary, we can take $\gamma \downarrow 0$ in \eqref{pr4}. As the left-hand side is free of $\gamma$, we get \eqref{geqn}. This completes the proof of Lemma \ref{slut}.

	\subsection{\bf Proof of Lemma \ref{tech}}\label{4sec4}
	We first show that $G_{\td}(0,0) \to g_{\zd}(0,0)$ as $n\to\infty$. We also derive the rate of convergence associated with it. 
	Towards this end, we define the boundary of a set $X$ in $\Z^d$ as $$\partial_{\zd} X := \{y\in X^c \mid y\mbox{ is a neighbour of } x\mbox{ for some } x\in X\}.$$
	Similarly we also define the boundary of a set $Y$ in $\td$ and denote it by $\partial_{\td} Y$. 
	Let $V=[1,n-2]^d\subset \zd$ and define $U=\Pi_n(V)$. Note that $U$ is properly contained in $\td$ in the sense that the boundary of $\widehat{U}$ is contained in $[0,n-1]^d \cap \zd$. This is important as it ensures that for all $x,y\in \td$, we have $g_{\td}^{U}(x,y)=g_{\zd}^{V}(\hat{x},\hat{y})$. See the discussion in Remark 1.8 of \cite{angelo} for more details.
	Hence
	\begin{equation}
	\begin{aligned}\label{bish}
	G_{\td}(0,0) & = G_{\td}\left(\lfloor\tfrac{n}2\rfloor,\lfloor\tfrac{n}2\rfloor\right) \\ & \stackrel{\ref{thr}}{=} g_{\td}^{U}(\lfloor\tfrac{n}2\rfloor,\lfloor\tfrac{n}2\rfloor)+ E_{\lfloor\tfrac{n}2\rfloor}^{\td}[G_{\td}(X_{T_{U}},\lfloor\tfrac{n}2\rfloor)] -\frac1{n^d}E_{\lfloor\tfrac{n}2\rfloor}^{\td}[T_{U}] \\ & = g_{\zd}^{V}(\lfloor\tfrac{n}2\rfloor,\lfloor\tfrac{n}2\rfloor)+ E_{\lfloor\tfrac{n}2\rfloor}^{\td}[G_{\td}(X_{T_{U}},\lfloor\tfrac{n}2\rfloor)] -\frac1{n^d}E_{\lfloor\tfrac{n}2\rfloor}^{\td}[T_{U}] \\ & \stackrel{\ref{one}}{=} g_{\zd}(\lfloor\tfrac{n}2\rfloor,\lfloor\tfrac{n}2\rfloor)-E_{\lfloor\tfrac{n}2\rfloor}^{\zd}[g_{\zd}(X_{T_V},\lfloor\tfrac{n}2\rfloor)\mathbb{I}_{\{T_V<\infty\}}]+\\ & \hspace{3cm}+E_{\lfloor\tfrac{n}2\rfloor}^{\td}[G_{\td}(X_{T_{U}},\lfloor\tfrac{n}2\rfloor)] -\frac1{n^d}E_{\lfloor\tfrac{n}2\rfloor}^{\td}[T_{U}].
	\end{aligned}
	\end{equation}
	We will now give an upper bound for $\ds E_{\lfloor\tfrac{n}2\rfloor}^{\td}[T_{U}]$. Since $U$ is properly contained in $\td$, we have $\ds E_{\lfloor\tfrac{n}2\rfloor}^{\td}[T_{U}] =E_{\lfloor\tfrac{n}2\rfloor}^{\zd}[T_{V}]$. Note that $V \subset  W:= \{x \in \zd \mid |x|< n\sqrt{d}  \}$. Then using estimates of $T_W$ (see equation (1.21) of \cite{lawler2013intersections}), we get that
	\begin{equation}\label{lawl}
	E_{\lfloor\tfrac{n}2\rfloor}^{\zd}[T_{V}] \le E_{\lfloor\tfrac{n}2\rfloor}^{\zd}[T_{W}] \le (n\sqrt{d}+2)^2 \le C_1n^2,
	\end{equation}
	for some constant $C_1$ depending on $d$ only. We further notice that $g_{\zd}(\lfloor\tfrac{n}2\rfloor,\lfloor\tfrac{n}2\rfloor)=g_{\zd}(0,0)$. Rearranging \eqref{bish}, and observing that $g_{\zd}$ and $\ds E_{\lfloor\tfrac{n}2\rfloor}^{\td}[T_{U}]$ are non-negative, we get that
	\begin{equation}\label{est}
	\begin{aligned}
	|G_{\td}(0,0)-g_{\zd}(0,0)|& \le E_{\lfloor\tfrac{n}2\rfloor}^{\zd}[g_{\zd}(X_{T_V},\lfloor\tfrac{n}2\rfloor)\mathbb{I}_{\{T_V<\infty\}}]+\\ & \hspace{2cm}+\left|E_{\lfloor\tfrac{n}2\rfloor}^{\td}[G_{\td}(X_{T_{U}},\lfloor\tfrac{n}2\rfloor)]\right| +\frac1{n^d}E_{\lfloor\tfrac{n}2\rfloor}^{\td}[T_{U}] \\ &  \stackrel{\eqref{lawl}}{\le} \sup_{z\in \partial_{\zd}V} g_{\zd}(z,\lfloor\tfrac{n}2\rfloor) + \sup_{z\in \partial_{\td}U} |G_{\td}(z,\lfloor\tfrac{n}2\rfloor)|+C_1n^{2-d} \\ & \le C_3n^{2-d}+C_2(\log(n))^{\frac{3d}{2}}n^{2-d}+C_1n^{2-d}\\ &  =\mathcal{O}((\log(n))^{\frac{3d}{2}}n^{2-d}).
	\end{aligned}
	\end{equation}
	The last inequality follows from the bounds on the Green's functions given in Lemma \ref{two} and Lemma \ref{four}. For simplicity we denote $v_n:=G_{\td}(0,0)$ and $v:=g_{\zd}(0,0)$. We now turn to the proof of Lemma \ref{tech}. Observe that 
	\begin{equation}\label{est1}
	|v_n-v|u_N(\delta)^2 = \mathcal{O}(\log(n)^{\frac{3d}2+1}n^{2-d})=o(1)
	\end{equation}
	where $u_N(\delta)=a_N\delta+b_N=\mathcal{O}(\sqrt{\log(n)})$ with $a_N$ and $b_N$ as defined in \eqref{con}.
	We note that by using the union bound and the standard inequality
	$$P(N(0,1)>a) \le \frac1{a\sqrt{2\pi}}e^{-a^2/2}$$ one can show that with high probability $\mathbf{A}_n$ and $\mathbf{B}_n$ occurs. To see this, observe that
	$$\begin{aligned}
	\mathbb{Q}_n(\mathbf{A}_n^{\mathsf{c}}) & \le \sum_{\alpha\in V_n\setminus R_n} \mathbb{P}^{\td}(\psi(\alpha)>\delta) = \sum_{\alpha\in V_n\setminus R_n} \mathbb{P}^{\td}(\Psi_{\td}(\Pi_n(\alpha))>u_N(\delta)) \\ & \le |V_n\setminus R_n|\exp\left(-\frac{u_N(\delta)^2}{2v_n}\right)\cdot \frac{\sqrt{v_n}}{\sqrt{2\pi}u_N(\delta)}
	\\ & = \left(n^d-(n-2n^{3/4})^d\right)\left[e^{-\frac{u_N(\delta)^2}{2v}} \frac{\sqrt{v}}{\sqrt{2\pi}u_N(\delta)}\right]\cdot e^{\frac{u_N(\delta)^2}{2v}-\frac{u_N(\delta)^2}{2v_n}} \frac{\sqrt{v_n}}{\sqrt{v}} \\ & \le C\frac{n^d-(n-2n^{3/4})^d}{n^d} \cdot \exp\left(-\frac{(v-v_n)u_N(\delta)^2}{2v_nv}\right)\sqrt{\frac{v_n}{v}} \\ & = o(1).
	\end{aligned}
	$$
	The last line in the above equation follows from \eqref{est1}. Thus $\ds\lim_{n\to\infty} \mathbb{Q}_n(\mathbf{A}_n^{\mathsf{c}})=0$. Similarly by applying the union bound again and the tail estimates of the Gaussian distribution, we get
	$$\begin{aligned}
	\mathbb{Q}_n(\mathbf{B}_n^{\mathsf{c}}) & \le \sum_{\alpha\in V_n\setminus R_n} \mathbb{P}^{\zd}(\phi(\alpha)>\delta) = \sum_{\alpha\in V_n\setminus R_n} \mathbb{P}^{\zd}\left(\varphi_{\zd}(\alpha)>u_N(\delta)\right)  \\ & \le |V_n\setminus R_n|\exp\left(-\frac{u_N(\delta)^2}{2v}\right)\cdot \frac{\sqrt{v}}{\sqrt{2\pi}u_N(\delta)}
	\\ & = |n^d-(n-2n^{3/4})^d|e^{-\frac{u_N(\delta)^2}{2v}} \frac{\sqrt{v}}{\sqrt{2\pi}u_N(\delta)} \\ & \le C\frac{n^d-(n-2n^{3/4})^d}{n^d} = o(1).
	\end{aligned}
	$$
	This shows that $\ds\lim_{n\to\infty} \mathbb{Q}_n(\mathbf{B}_n^{\mathsf{c}})=0$. For the event $\mathbf{C}_n(\gamma)$, we need to apply the coupling result of Lemma \ref{angl}. Observe that $b_N \le \sqrt{2v\log N}=\sqrt{2vd\log n}$. Thus $\ds a_N^2 = v^2b_N^{-2} \ge \frac{v}{2d\log (n)}$. Using this estimate, we have that
	$$\begin{aligned}
	\mathbb{Q}_n[(\mathbf{C}_n(\gamma))^{\mathsf{c}}] & \le \mathbb{Q}_n\left[\sup_{\alpha\in R_n} |\psi(\alpha)-\phi(\alpha)|\ge \gamma \right] \\ & \le \mathbb{Q}_n\left[\sup_{\alpha\in R_n}\left|\Psi_{\td}(\Pi_n(\alpha))-\varphi_{\zd}(\alpha)\right|\ge a_N\gamma\right]\\  & \stackrel{\ref{angl}}{\le} 4n^d\exp(-c_1a_N^2\gamma^2n^{c_2})\\ & \le 4n^d\exp\left(-c_1v\gamma^2\frac{n^{c_2}}{2d\log(n)}\right) = o(1).
	\end{aligned}
	$$
	Hence $\ds\lim_{n\to\infty} \mathbb{Q}_n[(\mathbf{C}_n(\gamma))^{\mathsf{c}}]=0$ which completes the proof of Lemma \ref{tech}.
	
	\section*{Acknowledgment} 
	
	We would like to thank Angelo Ab\"acherli, Alessandra Cipriani, and Parthanil Roy for helpful suggestions and corrections. The second author acknowledges the Matrics grant from the Department of Science and Technology for their support.

	\bibliographystyle{abbrvnat}			
	\bibliography{disbib}
\end{document}